\documentclass[a4paper,twoside,12pt,leqno]{article}
\usepackage{amsmath,amsthm,amssymb,enumerate}
\usepackage{eucal}
\usepackage{epsfig, float, afterpage, array}
\usepackage{hhline}
\usepackage{graphpap}

\swapnumbers
\theoremstyle{plain}
\newtheorem{theorem}{Theorem}[section]

\newtheorem{proposition}[theorem]{Proposition}

\theoremstyle{definition}
\newtheorem{definition}[theorem]{Definition}

\newtheorem{example}[theorem]{Example}

\numberwithin{equation}{theorem}
\numberwithin{figure}{theorem}

\def \Aut{\text{Aut}}
\def \Out{\text{Out}}
\def \M {\mathcal{M}}
\def \AM {\mathcal{AM}}

\def \ord{\text{ord}}
\def \A{\mathbb{A}}
\def \ordA{\text{ord}\mathbb{A}}
\def \SA{\text{S}\mathbb{A}}
\def \ordSA{\text{ordS}\mathbb{A}}
\def \ordD{\text{ord}D}
\def \Sym{\text{Sym}}

\def \reals {\mathbb{R}}

\begin{document}

\title{A combinatorial algorithm to compute presentations of mapping-class groups of orientable surfaces with one boundary component}

\author
{Llu\'{\i}s Bacardit\footnote{The research was
funded by MIC (Spain) through
Project MTM2008-01550.} }

\date{}

\maketitle


 \bigskip

\begin{abstract} 
We give an algorithm which computes a presentation for a subgroup, denoted $\AM_{g,1,p}$, of the automorphism group of a free group. It is known that $\AM_{g,1,p}$ is 
isomorphic to the mapping-class group of an orientable genus-$g$ surface with one boundary component and $p$ punctures. We 
define a variation of Auter space. 
\bigskip

{\footnotesize
\noindent \emph{{\normalfont 2010}\,Mathematics Subject Classification.} Primary: 57N05, 20F05;
Secondary: 20F28, 20F34.

\noindent \emph{Key words.} Mapping-class groups, presentations, automorphism groups, Auter space.}
\end{abstract}

\section{Introduction}

Let $S$ be an orientable genus-$g$ surface with $b$ boundary components and $p$ punctures. We denote by 
$\M(S)$ the group of isotopy classes of orientation-preserving homeomorphisms of $S$ which permute the set of punctures and pointwise fix the boundary components. Since the group $\M(S)$ only depends, up to isomorphism, on the genus $g$ of $S$, the 
number $b$ of boundary components of $S$ and the number $p$ of punctures of $S$, we denote $\M(S)$ by $\M_{g,b,p}$. We call 
$\M_{g,b,p}$ the {\it mapping-class group} of $S$.

Presentations for $\M_{g,b,p}$ were obtained after a sequence of papers started by Hatcher and Thurston~\cite{ht}, and followed by 
Harer~\cite{h}, Wajnryb~\cite{w1},\cite{w2}; Matsumoto~\cite{m} and, Labru\`{e}re and Paris~\cite{lp}. For $p=0$, Gervais~\cite{g} used \cite{w1} to deduce another 
presentations for $\M_{g,b,0}$. Before \cite{ht} very little was known about the presentation of $\M_{g,b,p}$. Birman and 
Hilden~\cite{bh} gave a presentation for $\M_{2,0,0}$, and McCool~\cite{mc} proved that $\M_{g,b,p}$ is finitely 
presented. 

Benvenuti~\cite{be} uses a variation of the curve complex, called ordered curve complex, to obtain presentations for 
$\M_{g,b,p}$ from an inductive process. This inductive process starts from presentation for the sphere and the torus with ``few'' boundary 
components and/or punctures. Hirose~\cite{hi} uses the curve complex and induction on $g$ and $b$ to deduce Gervais 
presentation. Both of these papers are independent of \cite{ht}.

Our algorithm is independent of \cite{ht}. We feel that our point of view goes back to McCool~\cite{mc}.
Section~\ref{sec:alg} contains the presentation given by our algorithm. This presentation has generators $ze_i,\ ze_ie_j$ where $z$ ranges over a finite set $\mathcal L$ and $e_i,\, e_j$ range over $z$. There are three type of relations:
\begin{enumerate}[(a).]
\item $ze_i=1,\,ze_{i'}e_{j'}=1$, for some generators $ze_i,\,ze_{i'}e_{j'}$.

\item $z_1e_ie_j = z_2e_{i'}e_{j'}$, for some generators $z_1e_ie_j,\, z_2e_{i'}e_{j'}$.

\item $ze_i\cdot ze_ie_j = ze_j\cdot ze_je_i$, for every generator $ze_ie_j$.
\end{enumerate}

\medskip

Armstrong, Forrest and Vogtmann~\cite{afv} give a new presentation for $\Aut(F_n)$, the automorphism group of 
the free group of rank $n$. This presentation for $\Aut(F_n)$ is obtained by studying the action of $\Aut(F_n)$ on a subcomplex of the spine of Auter space. Following Armstrong, Forrest and 
Vogtmann~\cite{afv}, we obtain our algorithm by studying the action of an algebraic analogous of $\M_{g,1,p}$ on a subcomplex of the spine of a variation of Auter space.

\section{Preliminaries}

Throughout the paper $n$ will be a non-negative integer, $F_n$ will be a free group of rank $n$, $\Aut(F_n)$ will be the automorphism group of $F_n$ and 
$\Out(F_n)$ will be the automorphism group of $F_n$ modulo inner automorphisms. Given $v,w\in F_n$, we denote by $[v,w]$ the element $v^{-1}w^{-1}vw$ 
of $F_n$. We denote by $[w]$ the conjugacy class of $w$.

\medskip

Let $S$ be an orientable genus-$g$ surface with $b$ boundary components and $p$ punctures. A homeomorphism $f$ of $S$ 
which fixes the basepoint of $\pi_1(S)$ and permutes the set of punctures of $S$ induces an automorphism $f_*\in \Aut(\pi_1(S))$. The isotopy class of $f$
defines an automorphism of $\pi_1(S)$ up to inner automorphisms, that is, an element of $\Out(\pi_1(S))$. For $(b,p)=(0,0)$, 
by Dehn-Nielsen-Baer Theorem, $\M_{g,0,0}$ is isomorphic to a index $2$ subgroup of $\Out(\pi_1(S))$. For $(g,p)\neq (0,0)$ by a modification of Dehn-Nielsen-Baer 
Theorem, $\M_{g,b,p}$ is isomorphic to an infinite index subgroup of $\Out(\pi_1(S))$.

Suppose now $b= 1$, that is, $S$ has exactly one boundary component. If we choose the 
basepoint of $\pi_1(S)$ to be a boundary point of $S$ and we restrict ourselves to homeomorphisms of $S$ which pointwise fix the boundary, then the isotopy class of a homeomorphism of $S$ defines an element of $\Aut(\pi_1(S))$. Since $S$ has one boundary component, the fundamental group of $S$ is 
a free group. We denote by $$\Sigma_{g,1,p}=\langle x_1,y_1,x_2,y_2,\ldots,x_g,y_g,t_1,t_2,\ldots, t_p \mid \,\,\,\rangle$$ 
a presentation of $\pi_1(S,*)$ where $*$ is a boundary point of $S$, for every $1\leq k\leq p$ the generator $t_k$ 
represents a loop around the $k$-th puncture of $S$ and the word $[x_1,y_1][x_2,y_2]\cdots[x_g,y_g]t_1t_2\cdots t_p$ 
represents a loop around the boundary component of $S$.

\begin{definition}
We denote by $\AM_{g,1,p}$ the subgroup of $\Aut(\Sigma_{g,1,p})$ consisting of automorphisms of $\Sigma_{g,1,p}$ which 
fix the word $[x_1,y_1][x_2,y_2]\cdots[x_g,y_g]t_1t_2\cdots t_p$ of $\Sigma_{g,1,p}$ and fix the set of conjugacy classes 
$[t_1^{-1}],[t_2^{-1}],\ldots,[t_p^{-1}]$ of $\Sigma_{g,1,p}$.
\qed
\end{definition}

Using a modification of Dehn-Nielsen-Baer Theorem, it can 
be proved that $\M_{g,1,p}$ is isomorphic to $\AM_{g,1,p}$, see~\cite{df} with some changes of notation and some 
different conventions. We call $\AM_{g,1,p}$ the {\it algebraic mapping-class group} of 
an orientable genus-$g$ surface with one boundary component and $p$ punctures.

\section{Auter space $\A_n$}

\begin{definition}\label{def:A_n}
Let $(\Gamma, v_0,\phi)$ be a $3$-tuple such that
\begin{enumerate}
\item\label{def:A_n:c1} $\Gamma$ is a finite connected graph with no separating edges.
\item\label{def:A_n:c2} $\Gamma$ is a metric graph with total volume $1$.
\item\label{def:A_n:c3} $v_0$ is a distinguished vertex of $\Gamma$.
\item\label{def:A_n:c4} Every vertex of $\Gamma$ but $v_0$ has valence at least $3$; $v_0$ has valence at least $2$.
\item\label{def:A_n:c5} $\phi:\pi_1(\Gamma,v_0) \to F_n$ is an isomorphism called ``marking''.
\end{enumerate}

A point in $\A_n$ is an equivalence class of $3$-tuples $(\Gamma,v_0,\phi)$, where $(\Gamma,v_0,\phi)$ is equivalent  to 
$(\Gamma',v_0',\phi')$ if there exists an isometry $h:\Gamma \to \Gamma$ such that $h(v_0)=v_0'$ and the isomorphism 
$h_*:\pi_1(\Gamma,v_0) \to \pi_1(\Gamma',v_0')$ satisfies $\phi = \phi'\circ h_*$.
\medskip

We call $\A_n$ Auter space.
\qed
\end{definition}

Auter space $\A_n$ was introduced by Hatcher and Vogtmann~\cite{hv} as an analogous for $\Aut(F_n)$ of Outer space. Often in the literature the marking is defined as $\phi^{-1}:F_n\to \pi_1(\Gamma,v_0)$.

\medskip

If $\Gamma$ has $k+1$ edges, then $(\Gamma,v_0,\phi)$ defines an open $k$-simplex of $\A_n$ denoted 
$\sigma(\Gamma,v_0,\phi)$. We can obtain $\sigma(\Gamma,v_0,\phi)$ by varying the length of the edges of $\Gamma$. The 
$k$-simplex $\sigma(\Gamma,v_0,\phi)$ can be parametrized by $\Delta^k$, the standard open $k$-simplex of $\reals^k$, 
as follows: $(\Gamma_s,v_0,\phi) \in \sigma(\Gamma,v_0,\phi)$ is the point of $\A_n$ such that the length of the edges of 
$\Gamma_s$ equal the barycentric coordinates of $s\in \Delta^k$. It is important that $\Delta^k$ is open. Since a non-trivial 
isometry of $\Gamma$ permutes same edges of $\Gamma$, such an isometry gives a non-trivial element of 
$H_1(\Gamma)$. Hence, every $s\in \Delta$ 
defines a different point of $\sigma(\Gamma,v_0,\phi)$.

Some faces of $\sigma(\Gamma,v_0,\phi)$ belong to $\A_n$. If an edge of $\Gamma$ is incident to two different vertices, then 
we can reduce the length of that edge to zero, and increase the length of the other edges, to obtain a new graph 
$\Gamma'$ with one edge minus. We say that we have 
{\it collapsed} one edge of $\Gamma$. We have a quotient map $\Gamma \to \Gamma'$ which defines a point 
$(\Gamma',v_0',\phi')$ of $\A_n$. We say that $\sigma(\Gamma',v_0',\phi')$ is a face of $\sigma(\Gamma,v_0,\phi)$. Faces of 
$\sigma(\Gamma',v_0',\phi')$ are faces of $\sigma(\Gamma,v_0,\phi)$. We cannot collapse an edge which is incident to a unique 
vertex. Hence, some face of $\sigma(\Gamma,v_0,\phi)$ are missing. In particular, $\A_n$ is not a simplicial complex.

There exists a deformation retract, denoted $\SA_n$, of $\A_n$ which is a simplicial complex. We can define $\SA_n$ as follows. For 
every simplex of $\A_n$, there exists a vertex of $\SA_n$. Two vertices of $\SA_n$ expand an edge if the simplex of $\A_n$ 
which defines one of the two vertices of $\SA_n$ is a face of the simplex of $\A_n$ which defines the other vertex of $\SA_n$; 
$i+1$ vertices of $\SA_n$ expand a $i$-simplex of $\SA_n$ if every pair of vertices expand an edge.

There exists a natural inclusion of $\SA_n$ into $\A_n$ by sending every vertex of $\SA_n$ to the barycenter of the 
corresponding simplex, and every $i$-simplex of $\SA_n$ to the convex hull of the corresponding barycenters. This inclusion 
is a deformation retract. See~\cite{hv}.

\medskip

Collapsing an edge of $\Gamma$ has an inverse process which {\it splits} a vertex of $\Gamma$ into two new vertices, 
and, the two new 
vertices are joined by a new edge. Often in the literature splitting of a vertex is called blowing up an edge. If $\tilde \Gamma$ is obtained from $\Gamma$ by splitting a vertex, 
then we can identify, in a natural way, every edge of $\Gamma$ with an edge of $\tilde \Gamma$. 
Collapsing the only edge of $\tilde \Gamma$ which is not identified with an edge of $\Gamma$ we obtain $\Gamma$.
There exists a quotient map $\tilde \Gamma \to \Gamma$. If $\tilde \Gamma$ is obtained from $\Gamma$ by splitting a vertex different from $v_0$, then the quotient map 
$\tilde \Gamma \to \Gamma$ defines a point $(\tilde \Gamma, \tilde v_0,\tilde \phi)$ of $\A_n$ 
If $\tilde \Gamma$ is obtained from $\Gamma$ by splitting $v_0$, 
then the point $(\tilde \Gamma, \tilde v_0,\tilde \phi)$ of $\A_n$ depends of the election, between the two possibilities,
of the new distinguished vertex $\tilde v_0$. 
The simplex 
$\sigma(\Gamma,v_0,\phi)$ is a face of $\sigma(\tilde \Gamma,\tilde v_0,\tilde \phi)$.

\medskip

We give a combinatorial 
definition of the topological type of $\Gamma$, that is, $\Gamma$ when we forget its metric. When we forget the 
metric of $\Gamma$ we can see $(\Gamma,v_0,\phi)$ as a vertex of $\SA_n$, in  fact, the vertex of $\SA_n$ 
defined by the simplex $\sigma(\Gamma,v_0,\phi)$ of $\A_n$. We translate 
to our combinatorial definition the processes of collapsing an edges and splitting a vertex. Our combinatorial definition 
of the topological type of $\Gamma$ is different from the one given in \cite{cv}.

\begin{definition}
Let \begin{enumerate}
     \item $V(\Gamma)$ be the set of vertices of $\Gamma$.
     \item $E(\Gamma)$ be the set of edges of $\Gamma$.
     \item $\overline E(\Gamma) = \{\overline e\mid e\in E(\Gamma)\}$ be a set disjoint with $E(\Gamma)$.
    \end{enumerate}
We extend $\overline{\phantom{\cdot\cdot }}\, $ to an involution of $E(\Gamma)\cup \overline E(\Gamma)$. We fix an orientation of every edge of $\Gamma$. We say that $e\in E(\Gamma)$ starts at $v_1\in V(\Gamma)$ and finishes 
at $v_2\in V(\Gamma)$ if $e$ is incident to $v_1$ and $v_2$; and $e$ is oriented from $v_1$ to $v_2$. In this case we say 
that $\overline e$ starts at $v_2$ and finishes at $v_1$.

\medskip

Given $v\in V(\Gamma)$, we define the following subset of $E(\Gamma)\cup \overline E(\Gamma)$. 
$$v^*=\{a\in E(\Gamma) \cup \overline E(\Gamma) \mid a \text{ starts at } v\}.$$
We set $V^*(\Gamma) = \{v^* \mid v\in V(\Gamma)\}$.

The topological type of $\Gamma$ is completely determined by $(V(\Gamma),E(\Gamma),V^*(\Gamma))$. 
\qed
\end{definition}

Notice that $v^*$ is the star of $v\in V(\Gamma)$ and $V^*(\Gamma)$ is a partition of 
$E(\Gamma) \cup \overline E(\Gamma)$. Condition \ref{def:A_n:c1} of Definition~\ref{def:A_n} can be translated by saying that $E(\Gamma)$ is finite and, for every $v\in V(\Gamma)$, 
there exist $a,b\in v^*$ such that $a\neq b,\overline b$ and $\overline a,\overline b \notin v^*$. Condition \ref{def:A_n:c4} 
of Definition~\ref{def:A_n} can be translated by saying that for every $v\in V(\Gamma)-\{v_0\}$, $v^*$ has at least $3$
 elements; $v_0^*$ has at least $2$ elements.   

\begin{definition}\label{def:collapsar}
Let $e\in E(\Gamma)$ such that $e\in v_1^*,\, \overline e\in v_2^*$, where $v_1,v_2\in V(\Gamma)$ and $v_1 \neq v_2$. We can collapse 
$e$. When we collapse the edge $e$ we have a graph with topological type 
$$(V(\Gamma)\cup \{v\}-\{v_1,v_2\}, E(\Gamma)-\{e\}, V^*(\Gamma)\cup \{v^*\}-\{v_1^*,v_2^*\})$$
where $v\notin V(\Gamma)$ and $v^*=v_1^*\cup v_2^*-\{e,\overline e\}$.
\qed
\end{definition}

\begin{definition}\label{def:partir}
Let $v\in V(\Gamma)-\{v_0\}$ and $A,B$ a partition of $v^*$ such that both $A$ and $B$ have at least two elements, there exists 
$a\in A$ such that $\overline a \notin A$ and there exists $b\in B$ such that $\overline b\notin B$. When we split the 
vertex $v$ with respect to $A$ and $B$ we have a graph with topological type $$(V(\Gamma)\cup\{v_1,v_2\}-\{v\},E(\Gamma)\cup \{e\}, V^*(\Gamma)\cup\{v_1^*,v_2^*\}-\{v^*\}),$$ where $v_1,v_2\notin V(\Gamma),\, e\notin E(\Gamma),\, v_1^*=A\cup \{e\}$ and $v_2^*=B\cup\{\overline e\}$. To split $v_0$ we have to choose between $v_1$ or $v_2$ as the new distinguished vertex. Since the distinguished vertex can have valence two, the subset of $v^*_0$ corresponding to the new distinguished vertex may have only one element.
\qed
\end{definition}

\section{Ordered Auter space $\ordA_{g,p}$}

Our motivation for defining ordered Auter space is that when a graph 
is embedded into an orientable surface, the star of every vertex of the graph which is mapped to an interior point of the surface 
gets a cyclic order, and, the star of a vertex which is mapped to a boundary point of the surface gets a linear order. When we want to 
collapse an edge or to split a vertex we have to do it respecting the orders of the stars.

\begin{definition}\label{def:qtuple}
Let $(\Gamma,v_0,\phi,\ord)$ be a $4$-tuple where $(\Gamma, v_0, \phi)$ satisfies conditions in Definition~\ref{def:A_n}, $\ord$ is a linear order of $v_0^*$ and a cyclic order of $v^*$ for every $v\in V(\Gamma)-\{v_0\}$.

Suppose $V(\Gamma) = \{v_0,v_1,v_2,\ldots,v_q\}$ and 
\begin{equation}\label{eq:ord}
\begin{array}{lcl}
\ord(v_0^*) & = & (a_1^0,a_2^0,\ldots, a_{r_0}^0),\\
\noalign{\vskip+0.05cm}
\ord(v_1^*) & = & (a_1^1,a_2^1,\ldots, a_{r_1}^1),\\
\noalign{\vskip+0.05cm}
\ord(v_2^*) & = & (a_1^2,a_2^2,\ldots, a_{r_2}^2),\\
\noalign{\vskip+0.05cm}
& \vdots\\
\noalign{\vskip+0.05cm}
\ord(v_q^*) & = & (a_1^q,a_2^q,\ldots,a_{r_q}^q).
\end{array}
\end{equation}
For $i\neq 0$, since $\ord(v_i^*)$ is cyclically ordered, the subindices of $\ord(v_i^*)$ are modulo $r_i$

We consider the following element of $\pi_1(\Gamma,v_0)$ and the following conjugacy classes of $\pi_1(\Gamma,v_0)$.
\begin{equation}
\begin{array}{lcl}\label{eq:par_ord}
w_0 & = & b_1^0b_2^0\cdots b_{l_0}^0,\\
\noalign{\vskip+0.05cm}
[w_1] & = & [b_1^1b_2^1\cdots b_{l_1}^1],\\
\noalign{\vskip+0.05cm}
[w_2] & = & [b_1^2b_2^2\cdots b_{l_2}^2],\\
\noalign{\vskip+0.05cm}
& \vdots\\
\noalign{\vskip+0.05cm}
[w_p] & = & [b_1^pb_2^p\cdots b_{l_p}^p],
\end{array}
\end{equation}
where $b_1^0=a_1^0$, for every $1\leq i\leq p,\, 1\leq j\leq l_i$ the subsequence $(\overline b_j^i,b^i_{j+1})$ appears in \eqref{eq:ord}, $b_{l_0}^{0}=\overline a_{r_q}^{\,q}$, and every element of $E(\Gamma)\cup E(\Gamma)$ appears exactly once in \eqref{eq:par_ord}.

We denote by $w(\Gamma,v_0,\ord)$ the set $\{w_0,[w_1],[w_2],\ldots,[w_p]\}$.
\qed
\end{definition}

\begin{example}\label{ex:1}
Let $(\Gamma,v_0,\ord)$ be a $3$-tuple where $V(\Gamma)=\{v_0,v_1,v_2\},\, E(\Gamma) = \{e_1,e_2,e_3,e_4,e_5\}$ and 
$$\begin{array}{lcl}
\ord(v_0^*) & = & (e_1,e_2),\\
\noalign{\vskip+0.05cm}
\ord(v_1^*) & = & (\overline e_1, e_3,\overline e_3,e_4,e_5),\\
\noalign{\vskip+0.05cm}
\ord(v_2^*) & = & (\overline e_2,\overline e_5,\overline e_4).
\end{array}$$
Then $w(\Gamma,v_0,\ord)= \{w_0,[w_1],[w_2],[w_3]\}$ where
$$\begin{array}{lcl}
w_0 & = & e_1e_3e_4\overline e_2,\\
\noalign{\vskip+0.05cm}
[w_1] & = & [\overline e_1e_2\overline e_5],\\
\noalign{\vskip+0.05cm}
[w_2] & = & [\overline e_3],\\
\noalign{\vskip+0.05cm}
[w_3] & = & [\overline e_4e_5].
\end{array}$$
\qed
\end{example}

Notice that for every $v\in V(\Gamma)$, $\ord(v^*)$ is completely determined by $w(\Gamma,v_0,\ord)$.

\begin{definition}
Let $(\Gamma,v_0,\phi,\ord)$ be a $4$-tuple such that $w(\Gamma,v_0,\ord)$ has $p$ conjugacy classes. 

We denote $\frac{n-p}{2}$ by $g$. We will see that $n-p$ is even. Hence, $g$ is a non-negative integer.

We define $\ordA_{g,p}$ as the space of equivalence classes of $4$-tuples $(\Gamma,v_0,\phi,\ord)$ such that $\phi: \pi_1(\Gamma,v_0)\to \Sigma_{g,1,p}$, $w(\Gamma,v_0,\ord) = \{w_0,[w_1],[w_2],\ldots,[w_p]\}$ and 
$$\begin{array}{lcl}
\phi(w_0) & = & [x_1,y_1][x_2,y_2]\cdots [x_g,y_g]t_1t_2\cdots t_p,\\
\noalign{\vskip+0.05cm}
\{\phi([w_1]),\phi([w_2]),\ldots,\phi([w_p])\} & = & \{[t_1^{-1}],[t_2^{-1}],\ldots,[t_p^{-1}]\}.
\end{array}$$

The $4$-tuples $(\Gamma,v_0,\phi,\ord)$ and $(\Gamma',v_0',\phi',\ord')$ represent the same point of $\ordA_{g,p}$ if there exists an isometry $h:\Gamma \to \Gamma'$ such that $h(v_0)=v_0'$, the isomorphism $h_*:\pi_1(\Gamma,v_0) \to \pi_1(\Gamma',v_0')$ satisfies $\phi = \phi'\circ h_*$, and $h:\Gamma \to \Gamma'$ preserves the orders, that is, $\ord(v^*) = (a_1,a_2,\ldots,a_r)$ implies $\ord'(h(v)^*) = (h(a_1),h(a_2),\ldots,h(a_r))$ for every $v\in V(\Gamma)$.

\medskip

We call $\ordA_{g,p}$ {\it ordered Auter space}. 

\medskip

We define $\ordSA_{g,p}$ for $\ordA_{g,p}$ as we defined $\SA_{n}$ for $\A_n$. In particular, $\ordSA_{g,p}$ is a simplicial complex, and, $\ordSA_{g,p}$ is a deformation retract of $\ordA_{g,p}$.
\qed
\end{definition}

The following definitions are based on Definition~\ref{def:collapsar} and Definition~\ref{def:partir}, respectively.

\begin{definition}\label{def:collapsar_ord}
Let $e\in E(\Gamma)$. Suppose $e = a_{k_1}^i,\, \overline e = a_{k_2}^j$, where $i\neq j$ and $1\leq k_1\leq r_i,\, 1\leq k_2\leq r_j$. Since $i\neq j$, we can collapse $e$. We can suppose $j\neq 0$. To adapt Definition~\ref{def:collapsar} to $\ordSA_{g,p}$ we set
\begin{align*}
\ord(v^*) = & (a_{k-1}^i,a_{k_2}^i,\ldots , a_{k_1-1}^i,\\
&\,\,\,\, a_{k_2+1}^j, a_{k_2+2}^j,\ldots , a_{r_j}^j, a_1^j, a_2^j, \ldots, a_{k_2-1}^j,\\
&\,\,\,\, a_{k_1+1}^i, a_{k_1+2}^i, \ldots, a_{r_1}^i).
\end{align*}
\qed
\end{definition}

\begin{example}\label{ex:2}
Let $(\Gamma,v_0,\ord)$ be as in Example~\ref{ex:1}. When we collapse $e_1$ we obtain $(\Gamma',v_0',\ord')$
such that $v_0' = v_0,\, V(\Gamma') = \{v_0, v_2\},\, E(\Gamma') = \{e_2,e_3,e_4,e_5\}$ and
$$\begin{array}{lcl}
\ord'(v_0^*) & = & (e_3,\overline e_3,e_4,e_5,e_2),\\
\noalign{\vskip+0.05cm}
\ord'(v_2^*) & = & (\overline e_2,\overline e_5,\overline e_4).
\end{array}$$
We have $w(\Gamma',v_0',\ord')= \{w'_0,[w'_1],[w'_2],[w'_3]\}$ where
$$\begin{array}{lcl}
w'_0 & = & e_3e_4\overline e_2,\\
\noalign{\vskip+0.05cm}
[w'_1] & = & [e_2\overline e_5],\\
\noalign{\vskip+0.05cm}
[w'_2] & = & [\overline e_3],\\
\noalign{\vskip+0.05cm}
[w'_3] & = & [\overline e_4e_5].
\end{array}$$
\qed
\end{example}

Let $(\Gamma,v_0,\phi,\ord)$ be a vertex of $\ordSA_{g,p}$. When we collapse an edge of $\Gamma$ according to Definition~\ref{def:collapsar_ord} we obtain $(\Gamma',v_0',\phi',\ord')$. As it is see in Example~\ref{ex:2}, $w(\Gamma',v_0',\phi',\ord')$ has $p$ conjugacy classes. Hence, $(\Gamma',v_0',\phi',\ord')$ is a vertex of $\ordSA_{g,p}$.

\begin{definition}\label{def:partir_ord}
Let $v\in V(\Gamma)$. Let $A,B$ be a partition of $v^*$. Suppose $\ord(v^*) = (a_1,a_2,\ldots, a_r)$ and $A=(a_{k_1}, a_{k_1+1}, \ldots, a_{k_2})$, where $1\leq k_1<k_2\leq r$. We can split the vertex $v$ with respect to $A,B$. To adapt Definition~\ref{def:partir} to $\ordSA_{g,p}$ we set
\begin{align*}
\ord(v_1^*) = &\,\, (e,a_{k_1},a_{k_1+1},\ldots, a_{k_2}),\\
\ord(v_2^*) = &\,\, (a_1,a_2,\ldots, a_{k_1-1},\\
&\,\,\,\,\,\, \overline e, a_{k_2+1}, a_{k_2+2},\ldots a_r).
\end{align*}
\qed
\end{definition}

\begin{example}\label{ex:3}
Let $(\Gamma,v_0,\ord)$ be as in Example~\ref{ex:1}. When we split $v_1$ with respect to $\{\overline e_3,e_4\},\, \{\overline e_1,e_3,e_5\}$ we obtain $(\tilde \Gamma,\tilde v_0,\tilde \ord)$
such that $\tilde v_0 = v_0,\, V(\tilde \Gamma) = \{v_0, v_{1,1},v_{1,2},v_2\},\, E(\tilde \Gamma) = \{ e,e_1,e_2,e_3,e_4,e_5\}$ and
$$\begin{array}{lcl}
\tilde \ord(v_0^*) & = & (e_1,e_2),\\
\noalign{\vskip+0.05cm}
\tilde \ord(v_{1,1}^*) & = & (e,\overline e_3,e_4),\\
\noalign{\vskip+0.05cm}
\tilde \ord(v_{1,2}^*) & = & (\overline e_1,e_3,\overline e,e_5),\\
\noalign{\vskip+0.05cm}
\tilde \ord(v_2^*) & = & (\overline e_2,\overline e_5,\overline e_4).
\end{array}$$
We have $w(\tilde \Gamma,\tilde v_0,\tilde \ord)= \{\tilde w_0,[\tilde w_1],[\tilde w_2],[\tilde w_3]\}$ where
$$\begin{array}{lcl}
\tilde w_0 & = & e_1e_3e_4\overline e_2,\\
\noalign{\vskip+0.05cm}
[\tilde w_1] & = & [\overline e_1e_2\overline e_5],\\
\noalign{\vskip+0.05cm}
[\tilde w_2] & = & [\overline e_3\overline e],\\
\noalign{\vskip+0.05cm}
[\tilde w_3] & = & [\overline e_4ee_5].
\end{array}$$
\qed
\end{example}

Let $(\Gamma,v_0,\phi,\ord)$ be a vertex of $\ordSA_{g,p}$. When we split a vertex of $\Gamma$ according to Definition~\ref{def:partir_ord} we obtain $(\tilde \Gamma,\tilde v_0,\tilde \phi,\tilde \ord)$. As it is see in Example~\ref{ex:3}, $w(\tilde \Gamma,\tilde v_0,\tilde \phi,\tilde \ord)$ has $p$ conjugacy classes. Hence, $(\tilde \Gamma,\tilde v_0,\tilde \phi,\tilde \ord)$ is a vertex of $\ordSA_{g,p}$.

\medskip

Since a graph satisfying Definition~\ref{def:A_n} can have at most $3n-2$ edges, the dimension of $\A_n$ is $3n-3$. On the other hand, $\A_n$ is not a manifold. The dimension of $\ordA_{g,p}$ is $6g+3p-3$; $\ordA_n$ is a manifold. 

\medskip

By~\cite[PROPOSITION 2.1]{hv} $\A_n$ is contractible. Since $\SA_n$ is a deformation retract of $\A_n$, we see that $\SA_n$ is contractible.

\begin{proposition}
$\ordA_{g,p}$ is contractible.
\end{proposition}
Hatcher and Vogtmann proof \cite[PROPOSITION 2.1]{hv} using spheres complexes. It is not clear how to translate to the context of spheres complexes an ordered graph. On the other hand, the proof of Culler and Vogtmann~\cite{cv} that Outer space is contractible can be applied to $\ordA_{g,p}$: adding a basepoint is straightforward, all the geometric arguments in \cite{cv} can be applied to $\ordA_{g,p}$ respecting the orders as Definition~\ref{def:collapsar_ord} and Definition~\ref{def:partir_ord} and McCool~\cite{mc2},~\cite{df2} proved that $\AM_{g,1,p}$ is generated by Nielsen automorphism which ``respect'' the orders (recall that Nielsen automorphisms are a special case of Whitehead automorphisms).

\medskip

Recall $n=2g+p$. There exists a natural map $\ordA_{g,p} \to \A_n$ which ``forgets'' the ordering, that is, $(\Gamma,v_0,\phi,\ord) \mapsto (\Gamma,v_0,\phi)$. Recall that $\ord$ is completely determined by $w(\Gamma,v_0,\ord)$. Since $\phi:\pi_1(\Gamma,v_0) \to \Sigma_{g,1,p}$ is an isomorphism, we have 
$$\begin{array}{ll}
w(\Gamma,v_0,\ord) = & \{\phi^{-1}([x_1,y_1][x_2,y_2]\cdots[x_g,y_g]t_1t_2\cdots t_p), \\
\noalign{\vskip+0.2cm}
&\quad\quad [\phi^{-1}(t_1^{-1})],[\phi^{-1}(t_2^{-1})],\ldots, [\phi^{-1}(t_p^{-1})])\}
\end{array}$$
Hence, the natural map $\ordA_{g,p} \to \A_n$ is injective.

\medskip

We want to see that $n-p$ is even.

For $n=1$, we have $p=1$ and $n-p=0$ is even. We do induction on $n$.

By Definition~\ref{def:collapsar_ord} we can collapse a maximal subtree of $(\Gamma,v_0,\ord)$. Hence, we can suppose that $V(\Gamma)=\{v_0\}$. Put $\ord(v_0^*) = (a_1,a_2,\ldots a_{2n})$ and $w(\Gamma,v_0,\ord) = \{w_0,[w_1],[w_2],\ldots,[w_p]\}$. Let $1\leq j\leq 2n$ such that $a_j = \overline a_1$. Then $w_0= a_1u$ in reduced form, for some $u\in F_n$. Let $(\Gamma',v_0',\ord')$ be obtained from $(\Gamma,v_0,\ord)$ be deleting the edges $a_1,a_j$. We have $\ord'(v_0') = (a_2,a_3,\ldots,a_{j-1},a_{j+1},a_{j+2},\ldots,a_{2n})$. We put $n'=n-1$ the rank of $\pi_1(\Gamma',v_0')$ and $p'$ the number of conjugacy classes of $w(\Gamma',v_0',\ord')$. By induction hypothesis $n'-p'$ is even.

If $w_0=a_1u'a_ju''$ cyclically reduced, then $w(\Gamma',v_0',\ord') = \{u'',[u'],[w_1],$ $[w_2],\ldots, [w_p]\}$. Notice that $u'\neq 1,\, u''\neq 1$ because $a_1u'a_ju''$ is cyclically reduced. Hence, $p'=p+1$ and $n-p=(n'+1)-(p'-1)=n'-p'+2$ is even.

If there exists $1\leq k\leq p$ such that $[w_k]=[a_jw_k']$, then $w(\Gamma',v_0',\ord') = \{w_k'u',[w_1], [w_2], \ldots, [w_{k-1}],[w_{k+1}],\ldots,[w_p]\}$. Hence, $p'=p-1$ and $n-p=(n'+1)-(p'+1)=n'-p'$ is even.

\section{The Degree Theorem}

Recall $\pi_1(\Gamma,v_0) \simeq F_n$.

\medskip

We denote the valence of $v\in V(\Gamma)$ by $|v^*|$.

\begin{definition}
The degree of $(\Gamma,v_0)$ is $2n-|v_0^*|$. Equivalently, the degree of $(\Gamma,v_0)$ is $\sum_{v\in V(\Gamma)-\{v_0\}} (|v^*|-2)$.
\qed
\end{definition}

To see the equivalence of the two definitions see~\cite[p. 636]{hv}. 

\medskip

From the first definition of the degree of $(\Gamma,v_0)$ we see that when we collapse an edge of $\Gamma$ which is not incident with $v_0$ the degree is preserved, and, when we collapse an edge of $\Gamma$ which is incident with $v_0$ the degree decreases. Hence, graphs of degree at most $i$ expand a subcomplex $D_i$ of $\SA_n$. Hatcher and Vogtmann~\cite{hv} proof the following.

\begin{theorem}
$D_i$ is $i$-dimensional and $(i-1)$-connected.
\end{theorem}

In particular, $D_2$ is a simply-connected $2$-complex.

\medskip

We define $\ordD_i$ for $\ordSA_{g,p}$ as we define $D_i$ for $\SA_n$.

\medskip

All the arguments of Hatcher and Vogtmann to proof \cite[THEOREM 3.3]{hv} can be applied to $\ordA_{g,p}$. In particular, 
what they call ``canonical splitting'' and ``sliding in the $\epsilon$-cone'' are combinations of splitting vertices and 
collapsing edges. We have the following.

\begin{theorem}\label{thm:simpl_ord}
$\ordD_i$ is $i$-dimensional and $(i-1)$-connected.
\end{theorem}

In particular, $\ordD_2$ is a simply-connected $2$-complex.

\section{The action of $\AM_{g,1,p}$ on $\ordA_n$}

Recall that $\Aut(F_n)$ acts on $\A_n$ by ``changing'' the markings: for every $\varphi \in \Aut(F_n)$ we define 
$\varphi\cdot (\Gamma,v_0,\phi) = (\Gamma,v_0,\varphi\circ \phi)$. This action restricts to $\SA_n$ and to $D_2$. The 
stabilizer of a vertex of $\SA_n$ by this action is a finite group which permutes some edges and invert some edge 
orientations. The quotient complex $\Aut(F_n)\backslash \SA_n$ is finite. See \cite[Section 3]{afv},~\cite[Section 5]{hv}. 
Armstrong, Forrest and Vogtmann~\cite{afv} apply a result of Brown~\cite{b} to $\Aut(F_n)\backslash D_2$ to compute a new 
presentation of $\Aut(F_n)$. Following this argument, we want to obtain a presentation of $\AM_{g,1,p}$.

\medskip

We can define an action of $\AM_{g,1,p}$ on $\ordA_{g,p}$ by ``changing'' the markings: for every $\varphi \in \AM_{g,1,p}$ 
we define $\varphi\cdot (\Gamma,v_0,\phi,\ord) = (\Gamma,v_0,\varphi\circ \phi,\ord)$. This action restricts to 
$\ordSA_{g,p}$ and to $\ordD_2$. The stabilizer of a vertex of $\ordSA_{g,p}$ by this action is trivial and the quotient 
complex $\AM_{g,1,p}\backslash \ordSA_{g,p}$ is finite, but much bigger than $\Aut(F_n)\backslash \SA_n$. By 
Theorem~\ref{thm:simpl_ord}, $\ordD_2$ is simply-connected. Hence, $\AM_{g,1,p}$ is isomorphic to the fundamental group of 
$\AM_{g,1,p}\backslash \ordD_2$. In the next section we give an algorithm which computes a presentation of the fundamental 
group of $\AM_{g,1,p}\backslash \ordD_2$.

\section{The algorithm}\label{sec:alg}

Recall $n=2g+p$.

\medskip

Vertices of $\ordSA_{g,p}$ are represented by $4$-tuples $(\Gamma,v_0,\phi,\ord)$ such that $w(\Gamma,v_0,\ord)$ has $p$ conjugacy classes. Recall that $\varphi \in \AM_{g,1,p}$ acts on $\ordSA_{g,p}$ by ``changing'' the marking, that is, $\varphi\cdot (\Gamma,v_0,\phi,\ord) = (\Gamma,v_0,\varphi\circ \phi,\ord)$. Hence, the quotient map $\ordSA_{g,p} \to \AM_{g,1,p}\backslash \ordSA_{g,p},\, (\Gamma,v_0,\phi,\ord) \mapsto (\Gamma,v_0,\ord)$ ``forgets'' the marking.
We can represent vertices of $\AM_{g,1,p}\backslash \ordSA_{g,p}$ by $3$-tuples $(\Gamma,v_0,\ord)$ such that $w(\Gamma,v_0,\ord)$ has $p$ conjugacy classes. Vertices of the subcomplex $\AM_{g,1,p}\backslash \ordD_2$ can be represented by $3$-tuples $(\Gamma,v_0,\ord)$ such that $(\Gamma,v_0)$ has degree at most $2$.

We want to compute a presentation for the fundamental group of complex $\AM_{g,1,p}\backslash \ordD_2$. Recall that the degree of $(\Gamma,v_0)$ can be defined as $\sum_{v\in V(\Gamma)-\{v_0\}}(|v^*|-2)$. Hence, if $(\Gamma,v_0)$ has degree $2$ then $\Gamma$ has at most three vertices: $v_0$ and two more vertices of valence $3$.

Let $\mathcal L$ be a list of vertices $(\Gamma,v_0,\ord)$ of $\AM_{g,1,p}\backslash \ordD_2$ such that $\Gamma$ has $3$ vertices.

Let $z=(\Gamma,v_0,\ord)$ be an element of $\mathcal L$. Suppose $E(\Gamma) = \{e_1,e_2,\ldots e_k\}$.

We construct a tree $T(z)$ as follows. There exists a vertex $z$ of $T(z)$. Let $e_i$ be an edge of $\Gamma$ which can be collapsed, that is, $e_i$ is incident to two different vertices. When we collapse $e_i$ we have a quotient $3$-tuple $z^i=(\Gamma^i,v_0^i,\ord^i)$. There exists a vertex $z^i$ of $T(z)$ and an edge $ze_i$ of $T(z)$ from $z$ to $z^i$. We identify edges of $z^i$ with edges of $z$. Let $e_j$ be an edge of $\Gamma^i$ which can be collapsed. When we collapse $e_j$ in $\Gamma^i$ we have a quotient $3$-tuple $z^{(i,j)}=(\Gamma^{(i,j)},v_0^{(i,j)},\ord^{(i,j)})$. There exists a vertex $z^{(i,j)}$ of $T(z)$ and an edge $ze_ie_j$ from $z^i$ to $z^{(i,j)}$. We repeat this process for every edge which can be collapsed.

Our generating set for the fundamental group of $\AM_{g,1,p}\backslash \ordD_2$ is the set of edges of $T(z)$, where $z$ ranges over $\mathcal L$.

\medskip

The group $\Sym_k\times C_2^{\times k}$ acts on $E(\Gamma)\cup \overline E(\Gamma)$ by permuting edges ($C_2^{\times k}$ is the Cartesian product of $k$ copies of the cyclic group of order $2$). Hence, $\Sym_k\times C_2^{\times k}$ acts on the set of $3$-tuples $(\Gamma,v_0,\ord)$ by permuting edges and inverting edge orientations. Two $3$-tuples $(\Gamma,v_0,\ord)$ and $(\Gamma',v_0',\ord')$ represent the same vertex of $\AM_{g,1,p}\backslash \ordSA_{g,p}$ if and only if they are in the same orbit by the action of $\Sym_k\times C_2^{\times k}$. Since every vertex of $T(z)$ is a $3$-tuple $(\Gamma,v_0,\ord)$, we see that $\Sym_k\times C_2^{\times k}$ acts on $T(z)$. We can identify $(\Sym_k\times C_2^{\times k})\backslash T(z)$ with the $1$-skeleton of a subcomplex of $\AM_{g,1,p}\backslash \ordD_2$. We can identify $$(\Sym_k\times C_2^{\times k})\backslash \Big(\bigcup_{z\in \mathcal L}T(z)\Big)$$ with the $1$-skeleton of a subcomplex of $\AM_{g,1,p}\backslash \ordD_2$. 

\medskip

We attach some $2$-cells to $(\Sym_k\times C_2^{\times k})\backslash \big(\bigcup_{z\in \mathcal L}T(z)\big)$. If there exists the generator $ze_ie_j$, we attach a $2$-cell though the egde-path $ze_i, ze_ie_j, \overline{ze_je_i}$, $\overline{ze_j}$. With these attached $2$-cells, the $2$-complex $(\Sym_k\times C_2^{\times k})\backslash \big(\bigcup_{z\in \mathcal L}T(z)\big)$ is homeomorphic to $\AM_{g,1,p}\backslash \ordD_2$. We fix a maximal subtree of $(\Sym_k\times C_2^{\times k})\backslash \big(\bigcup_{z\in \mathcal L}T(z)\big)$. 

\medskip

Our presentation for the fundamental group of $\AM_{g,1,p}\backslash \ordD_2$ has three types of relations:
\begin{enumerate}[(a).]
\item $ze_i=1,\,ze_{i'}e_{j'}=1$ if the edges $ze_i,\,ze_{i'}e_{j'}$ are in our maximal subtree.

\item $z_1e_ie_j = z_2e_{i'}e_{j'}$ if the generator $z_1e_ie_j$ exists and $g\cdot z_2^{i'} = z_1^{i}$ for some $g \in \Sym_k\times C_2^{\times k}$ such that either $g\cdot e_{j'} = e_j$ or $g\cdot e_{j'} = \overline e_j$.

\item $ze_i\cdot ze_ie_j = ze_j\cdot ze_je_i$ if there exists the generator $ze_ie_j$.
\end{enumerate}

We illustrate the algorithm with two easy examples. The main difficulty of the algorithm is to find $\mathcal L$. Once $\mathcal L$ is known our, it is straightforward to apply the algorithm. Example~\ref{ex:alg_2} shows that the algorithm can be applied in ``pieces'', each piece corresponding to an element of $\mathcal L$.

\begin{example}
We take $(g,p) = (1,0)$. The list $\mathcal{L}$ has a single element. We can represent the element of $\mathcal L$ by $$z=(V(\Gamma),E(\Gamma),V^*(\Gamma),\ord)=(\{v_0,v_1,v_2\},\{e_1,e_2,e_3,e_4\},\{v_0^*,v_1^*, v_2^*\},\ord),$$ where $\ord(v_0^*)=(e_1,e_2),\, \ord(v_1^*)=(\overline e_1, e_3,e_4) \text{ and } \ord(v_2^*)=(\overline e_2,\overline e_3, \overline e_4)$. To simplify the notation we put $$z=\ord(v_0^*); \ord(v_1^*), \ord(v_2^*)= (e_1,e_2);(\overline e_1, e_3,e_4),(\overline e_2,\overline e_3, \overline e_4).$$ We can collapse all $4$ edges of $z$ and we have 
\begin{align*}
z^1 & =(e_3,e_4,e_2);(\overline e_2,\overline e_3, \overline e_4),\\
z^2 & =(e_1,\overline e_3, \overline e_4);(\overline e_1, e_3,e_4),\\
z^3 & =(e_1,e_2);(\overline e_1, \overline e_4, \overline e_2, e_4),\\
z^4 & =(e_1,e_2);(\overline e_1, e_3,\overline e_2,\overline e_3).
\end{align*}
We see $z^1=g^{2,1}\cdot z^2,\, z^3=g^{4,3}\cdot z^4$, where $g^{2,1},\, g^{4,3}\in \Sym_4\times C_2^{\times 4}$ and

$$g^{2,1}=\left\{
\begin{array}{ll}
e_1 & \mapsto e_3,\\
e_3 & \mapsto \overline e_4,\\
e_4 & \mapsto \overline e_2,
\end{array}\right.\quad \quad \text{and} \quad\quad 
g^{4,3}=\left\{
\begin{array}{ll}
e_1 & \mapsto e_1,\\
e_2 & \mapsto e_2,\\
e_3 & \mapsto \overline e_4.
\end{array}\right.$$

We can collapse some edges of $z^1$ and $z^3$ and we have  
$$\begin{array}{ll}
z^{(1,3)} & =(\overline e_4, \overline e_2,e_4,e_2),\\
z^{(1,4)} & =(e_3,\overline e_2,\overline e_3,e_2),\\
z^{(1,2)} & =(e_3,e_4,\overline e_3, \overline e_4);
\end{array}
\quad\quad \text{and} \quad\quad
\begin{array}{ll}
z^{(3,1)} & =(\overline e_4, \overline e_2, e_4,e_2),\\
z^{(3,2)} & =(e_1,e_4,\overline e_1, \overline e_4).
\end{array}
$$
We see $z^{(1,2)},\, z^{(1,3)},\, z^{(1,4)},\, z^{(3,1)}$ and $z^{(3,2)}$ are in the same orbit by $\Sym_4\times C_2^{\times 4}$. Hence, they represent the same vertex of $(\Sym_4\times C_2^{\times 4})\backslash T(z)$.

\medskip

We take the maximal subtree of $(\Sym_4\times C_2^{\times 4})\backslash T(z)$ with edges  $ze_1, ze_3$ and  $ze_1e_2$. Then $\AM_{1,0}$ has presentation with generators:
$$\begin{array}{l}
ze_1,\, ze_2,\, ze_3,\, ze_4,\\
ze_1e_3,\, ze_1e_4,\, ze_1e_2,\\
ze_2e_1,\, ze_2e_3,\, ze_2e_4,\\
ze_3e_1,\, ze_3e_2,\\
ze_4e_1,\, ze_4e_2;
\end{array}$$
and relations: 
$$
\begin{array}{l}
ze_1= 1,\, ze_3 = 1,\, ze_1e_2=1,\\
ze_2e_1=ze_1e_3,\, ze_2e_3=ze_1e_4,\, ze_2e_4=ze_1e_2,\, ze_4e_1=ze_3e_1,\, ze_4e_2=ze_3e_2,\\
ze_1\cdot ze_1e_2 = ze_2\cdot ze_2e_1,\, ze_1\cdot ze_1e_3 = ze_3\cdot ze_3e_1,\, ze_1\cdot ze_1e_4 = ze_4\cdot ze_4e_1,\\
\quad\quad ze_2\cdot ze_2e_3= ze_3\cdot ze_3e_2,\,ze_2\cdot ze_2e_4= ze_4\cdot ze_4e_2.
\end{array}
$$
An easy simplification shows $\AM_{1,1,0}=\langle ze_2,ze_4\mid ze_2\cdot ze_2 = ze_4\cdot ze_2\cdot ze_4 \rangle$.
\qed
\end{example}

\begin{example}\label{ex:alg_2}
We take $(g,p) = (0,3)$. The list $\mathcal{L}$ is 
\begin{align*}
z_1 &= (e_1,e_2,e_3,e_4);(\overline e_1,e_5,\overline e_2),(\overline e_3,\overline e_5,\overline e_4),\\
z_2 &= (e_1,e_2,e_3,e_4);(\overline e_1,\overline e_4,e_5),(\overline e_2,\overline e_5,\overline e_3),\\
z_3 &= (e_1,\overline e_1,e_2,e_3);(\overline e_2,e_4,e_5),(\overline e_3,\overline e_5,\overline e_4),\\
z_4 &= (e_1,e_2,\overline e_2,e_3);(\overline e_1,e_4,e_5),(\overline e_3,\overline e_5,\overline e_4),\\
z_5 &= (e_1,e_2,e_3,\overline e_3);(\overline e_1,e_4,e_5),(\overline e_2,\overline e_5,\overline e_4),\\
z_6 &= (e_1,e_2,e_3,\overline e_1);(\overline e_2,e_4,e_5),(\overline e_3,\overline e_5,\overline e_4).
\end{align*}

For $z_1$ we have
\begin{align*}
z_1 &= (e_1,e_2,e_3,e_4);(\overline e_1,e_5,\overline e_2),(\overline e_3,\overline e_5,\overline e_4),\\
& z_1^1 = (e_5,\overline e_2,e_2,e_3,e_4);(\overline e_3,\overline e_5,\overline e_4),\\
& z_1^2 = (e_1,\overline e_1,e_5,e_3,e_4);(\overline e_3,\overline e_5,\overline e_4),\\
& z_1^3 = (e_1,e_2,\overline e_5,\overline e_4,e_4);(\overline e_1,e_5,\overline e_2),\\
& z_1^4 = (e_1,e_2,e_3,\overline e_3,\overline e_5);(\overline e_1,e_5,\overline e_2),\\
& z_1^5 = (e_1,e_2,e_3,e_4);(\overline e_1,\overline e_4,\overline e_3,\overline e_2).
\end{align*}

Generators for $z_1$ are:
\begin{align*}
&z_1e_1,\, z_1e_2,\, z_1e_3,\, z_1e_4,\, z_1e_5,\\
&z_1e_1e_5,\, z_1e_1e_3,\, z_1e_1e_4,\\
&z_1e_2e_5,\, z_1e_2e_3,\, z_1e_2e_4,\\
&z_1e_3e_1,\, z_1e_3e_2,\, z_1e_3e_5,\\
&z_1e_4e_1,\, z_1e_4e_2,\, z_1e_4e_5,\\
&z_1e_5e_1,\, z_1e_5e_2,\, z_1e_5e_3,\, z_1e_5e_4.
\end{align*}

We have 
\begin{align*}
z_1^{(1,5)} &= (\overline e_4,\overline e_3,\overline e_2,e_2,e_3,e_4),\\
z_1^{(1,3)} &= (e_5,\overline e_2,e_2,\overline e_5,\overline e_4,e_4),\\
z_1^{(1,4)} &= (e_5,\overline e_2,e_2,e_3,\overline e_3,\overline e_5),\\
z_1^{(2,5)} &= (e_1,\overline e_1,\overline e_4,\overline e_3,e_3,e_4),\\
z_1^{(2,3)} &= (e_1,\overline e_1,e_5,\overline e_5,\overline e_4,e_4),\\
z_1^{(2,4)} &= (e_1,\overline e_1,e_5,e_3,\overline e_3,\overline e_5),\\
z_1^{(3,5)} &= (e_1,e_2,\overline e_2,\overline e_1,\overline e_4,e_4),\\
z_1^{(4,5)} &= (e_1,e_2,e_3,\overline e_3,\overline e_2,\overline e_1).
\end{align*}
We see that $z_1^{(2,4)} = g\cdot z_1^{(2,5)},\, z_1^{(3,5)} = g'\cdot z_1^{(1,3)},\, z_1^{(4,5)} = g''\cdot z_1^{(1,5)}$ for some $g,g',g''\in \Sym_5\times C_2^{\times 5}$.

\medskip

Relations for $z_1$ are:
\begin{align*}
&z_1e_1=1,\, z_1e_2=1,\, z_1e_3=1,\, z_1e_4=1,\, z_1e_5=1,\\
&z_1e_1e_5=1,\, z_1e_1e_3=1,\, z_1e_1e_4=1,\, z_1e_2e_5=1,\, z_1e_2e_3= 1,\\
&z_1e_1\cdot z_1e_1e_5 = z_1e_5\cdot z_1e_5e_1,\,  z_1e_1\cdot z_1e_1e_3 = z_1e_3\cdot z_1e_3e_1,\\
& \quad\quad z_1e_1\cdot z_1e_1e_4 = z_1e_4\cdot z_1e_4e_1,\\
&z_1e_2\cdot z_1e_2e_5 = z_1e_5\cdot z_1e_5e_2,\, z_1e_2\cdot z_1e_2e_3 = z_1e_3\cdot z_1e_3e_2,\\ 
& \quad\quad z_1e_2\cdot z_1e_2e_4 = z_1e_4\cdot z_1e_4e_2,\\
&z_1e_3\cdot z_1e_3e_5 = z_1e_5\cdot z_1e_5e_3,\\
&z_1e_4\cdot z_1e_4e_5 = z_1e_5\cdot z_1e_5e_4.
\end{align*}

An easy simplification shows that for $z_1$ generators are $z_1e_2e_4,\, z_1e_3e_5,\, z_1e_4e_5$ and for $z_1$ there are no relations.

\medskip

From $z_2$ we have

\begin{align*}
z_2 &= (e_1,e_2,e_3,e_4);(\overline e_1,\overline e_4,e_5),(\overline e_2,\overline e_5,\overline e_3),\\
& z_2^1 = (\overline e_4,e_5,e_2,e_3,e_4);(\overline e_2,\overline e_5,\overline e_3),\\
& z_2^2 = (e_1,\overline e_5,\overline e_3,e_3,e_4);(\overline e_1,\overline e_4,e_5),\\
& z_2^3 = (e_1,e_2,\overline e_2,\overline e_5,e_4);(\overline e_1,\overline e_4,e_5),\\
& z_2^4 = (e_1,e_2,e_3,e_5,\overline e_1);(\overline e_2,\overline e_5,\overline e_3),\\
& z_2^5 = (e_1,e_2,e_3,e_4);(\overline e_1,\overline e_4,\overline e_3,\overline e_2).
\end{align*}

Generators for $z_2$ are:
\begin{align*}
&z_2e_1,\, z_2e_2,\, z_2e_3,\, z_2e_4,\, z_2e_5,\\
&z_2e_1e_5,\, z_2e_1e_2,\, z_2e_1e_3,\\
&z_2e_2e_1,\, z_2e_2e_5,\, z_2e_2e_4,\\
&z_2e_3e_1,\, z_2e_3e_5,\, z_2e_3e_4,\\
&z_2e_4e_2,\, z_2e_4e_3,\, z_2e_4e_5,\\
&z_2e_5e_1,\, z_2e_5e_2,\, z_2e_5e_3,\, z_2e_5e_4.
\end{align*}

We see $z_1^4=g_{2,1}^{2,4}\cdot z_2^2, \, z_1^1 = g_{2,1}^{3,1}\cdot z_2^3, \, z_2^1 = g_{2,2}^{4,1}\cdot z_2^4, \, z_1^5 = g_{2,1}^{5,5}\cdot z_2^5$, where
$$
\begin{array}{lll}
g_{2,1}^{2,4}=\left\{\begin{array}{lll}
e_1 & \mapsto & e_1,\\
e_3 & \mapsto & \overline e_3,\\
e_4 & \mapsto & \overline e_5,\\
e_5 & \mapsto & \overline e_2;
\end{array}\right.
& 
g_{2,1}^{3,1}=\left\{\begin{array}{lll}
e_1 & \mapsto & e_5,\\
e_2 & \mapsto & \overline e_2,\\
e_4 & \mapsto & e_4,\\
e_5 & \mapsto & \overline e_3;
\end{array}\right.
&
g_{2,2}^{4,1} =\left\{\begin{array}{lll}
e_1 & \mapsto & \overline e_4,\\
e_2 & \mapsto & e_5,\\
e_3 & \mapsto & e_2,\\
e_5 & \mapsto & e_3;
\end{array}\right.
\\
\noalign{\vskip+0.25cm}
g_{2,1}^{5,5}=\left\{\begin{array}{lll}
e_1 & \mapsto & e_1,\\
e_2 & \mapsto & e_2,\\
e_3 & \mapsto & e_3,\\
e_4 & \mapsto & e_4.
\end{array}\right.
\end{array}
$$

Relations for $z_2$ are:
\begin{align*}
&z_2e_1=1,z_2e_2=1,\\
&z_2e_2e_1 = z_1e_4e_1,\, z_2e_2e_5 = z_1e_4e_2,\, z_2e_2e_4 = z_1e_4e_5,\\
&z_2e_3e_1 = z_1e_1e_5,\, z_2e_3e_5 = z_1e_1e_3,\, z_2e_3e_4 = z_1e_1e_4,\\
&z_2e_4e_2 = z_2e_1e_5,\, z_2e_4e_3 = z_2e_1e_2,\, z_2e_4e_5 = z_2e_1e_3,\\
&z_2e_5e_1 = z_1e_5e_1,\, z_2e_5e_2 = z_1e_5e_2,\, z_2e_5e_3 = z_1e_5e_3,\, z_2e_5e_4 = z_1e_5e_4,\\
&z_2e_1\cdot z_2e_1e_5 = z_2e_5\cdot z_2e_5e_1,\,  z_2e_1\cdot z_2e_1e_2 = z_2e_2\cdot z_2e_2e_1,\\
& \quad\quad z_2e_1\cdot z_2e_1e_3 = z_2e_3\cdot z_2e_3e_1,\\
&z_2e_2\cdot z_2e_2e_5 = z_2e_5\cdot z_2e_5e_2,\, z_2e_2\cdot z_2e_2e_4 = z_2e_4\cdot z_2e_4e_2,\\ 
&z_2e_3\cdot z_2e_3e_5 = z_2e_5\cdot z_2e_5e_3,\, z_2e_3\cdot z_2e_3e_4 = z_2e_4\cdot z_2e_4e_3,\\
&z_2e_4\cdot z_2e_4e_5 = z_2e_5\cdot z_2e_5e_4.
\end{align*}

An easy simplification shows no new generators are needed, the relations $z_1e_4e_5=z_1e_2e_4\cdot z_1e_3e_5\cdot z_1e_2e_4,\, z_1e_4e_5\cdot z_1e_3e_5=z_1e_2e_4\cdot z_1e_4e_5$ are needed.

\medskip

From $z_3$ we have

\begin{align*}
z_3 &= (e_1,\overline e_1,e_2,e_3);(\overline e_2,e_4,e_5),(\overline e_3,\overline e_5,\overline e_4),\\
& z_3^2 = (e_1,\overline e_1,e_4,e_5,e_3);(\overline e_3,\overline e_5,\overline e_4),\\
& z_3^3 = (e_1,\overline e_1,e_2,\overline e_5,\overline e_4);(\overline e_2,e_4,e_5),\\
& z_3^4 = (e_1,\overline e_1,e_2,e_3);(\overline e_2,\overline e_3,\overline e_5,e_5),\\
& z_3^5 = (e_1,\overline e_1,e_2,e_3);(\overline e_2,e_4,\overline e_4,\overline e_3).
\end{align*}

Generators for $z_3$ are:
\begin{align*}
&z_3e_2,\, z_3e_3,\, z_3e_4,\, z_3e_5,\\
&z_3e_2e_4,\, z_3e_2e_5,\, z_3e_2e_3,\\
&z_3e_3e_2,\, z_3e_3e_5,\, z_3e_3e_4,\\
&z_3e_4e_2,\, z_3e_4e_3,\\
&z_3e_5e_2,\, z_3e_5e_3. \\
\end{align*}

We see $z_1^2=g_{3,1}^{2,2}\cdot z_3^2, \, z_1^2 = g_{3,1}^{3,2}\cdot z_3^3$, where
$$
\begin{array}{ll}
g_{3,1}^{2,2}=\left\{\begin{array}{lll}
e_1 & \mapsto & e_1,\\
e_3 & \mapsto & e_4,\\
e_4 & \mapsto & e_5,\\
e_5 & \mapsto & e_3;\\
\end{array}\right.
&
g_{3,1}^{3,2}=\left\{\begin{array}{lll}
e_1 & \mapsto & e_1,\\
e_2 & \mapsto & e_5,\\
e_4 & \mapsto & \overline e_4,\\
e_5 & \mapsto & \overline e_3.
\end{array}\right.
\end{array}
$$

Relations for $z_3$ are:
\begin{align*}
&z_3e_2,\, z_3e_4=1,\, z_3e_5=1,\\
&z_3e_2e_4 = z_1e_2e_5,\, z_3e_2e_5 = z_1e_2e_3,\,z_3e_2e_3 = z_1e_2e_4,\\
&z_3e_3e_2 = z_1e_2e_5,\, z_3e_3e_5 = z_1e_2e_3,\,z_3e_3e_4 = z_1e_2e_4,\\
&z_3e_2\cdot z_3e_2e_4 = z_3e_4\cdot z_3e_4e_2,\,  z_3e_2\cdot z_3e_2e_5 = z_3e_5\cdot z_3e_5e_2,\\
& \quad\quad z_3e_2\cdot z_3e_2e_3 = z_3e_3\cdot z_3e_3e_2,\\
&z_3e_3\cdot z_3e_3e_5 = z_3e_5\cdot z_3e_5e_3,\, z_3e_3\cdot z_3e_3e_4 = z_3e_4\cdot z_3e_4e_3.
\end{align*}

An easy simplification shows that neither new generators nor new relations are needed.

\medskip

From $z_4$ we have

\begin{align*}
z_4 &= (e_1,e_2,\overline e_2,e_3);(\overline e_1,e_4,e_5),(\overline e_3,\overline e_5,\overline e_4),\\
& z_4^1 = (e_4,e_5,e_2,\overline e_2,e_3);(\overline e_3,\overline e_5,\overline e_4),\\
& z_4^3 = (e_1,e_2,\overline e_2,\overline e_5,\overline e_4);(\overline e_1,e_4,e_5),\\
& z_4^4 = (e_1,e_2,\overline e_2,e_3);(\overline e_1,\overline e_3,\overline e_5,e_5),\\
& z_4^5 = (e_1,e_2,\overline e_2,e_3);(\overline e_1,e_4,\overline e_4,\overline e_3).
\end{align*}

Generators for $z_4$ are:
\begin{align*}
&z_4e_1,\, z_4e_3,\, z_4e_4,\, z_4e_5,\\
&z_4e_1e_4,\, z_4e_1e_5,\, z_4e_1e_3,\\
&z_4e_3e_1,\, z_4e_3e_5,\, z_4e_3e_4,\\
&z_4e_4e_1,\, z_4e_4e_3,\\
&z_4e_5e_1,\, z_4e_5e_3. \\
\end{align*}

We see $z_1^4=g_{4,1}^{1,4}\cdot z_4^1, \, z_1^1 = g_{4,1}^{3,1}\cdot z_4^3$, where
$$
\begin{array}{ll}
g_{4,1}^{1,4}=\left\{\begin{array}{lll}
e_2 & \mapsto & e_3,\\
e_3 & \mapsto & \overline e_5,\\
e_4 & \mapsto & e_1,\\
e_5 & \mapsto & e_2;\\
\end{array}\right.
& 
g_{4,1}^{3,1}=\left\{\begin{array}{lll}
e_1 & \mapsto & e_5,\\
e_2 & \mapsto & e_2,\\
e_4 & \mapsto & \overline e_4,\\
e_5 & \mapsto & \overline e_3.
\end{array}\right.
\end{array}
$$

Relations for $z_4$ are:
\begin{align*}
&z_4e_4=1,\, z_4e_4=1,\, z_4e_5=1,\\
&z_4e_1e_4 = z_1e_4e_1,\, z_4e_1e_5 = z_1e_4e_2,\,z_4e_1e_3 = z_1e_4e_5,\\
&z_4e_3e_1 = z_1e_1e_5,\, z_4e_3e_5 = z_1e_1e_3,\,z_4e_3e_4 = z_1e_1e_4,\\
&z_4e_1\cdot z_4e_1e_4 = z_4e_4\cdot z_4e_4e_1,\,  z_4e_1\cdot z_4e_1e_5 = z_4e_5\cdot z_4e_5e_1,\\
& \quad\quad z_4e_1\cdot z_4e_1e_3 = z_4e_3\cdot z_4e_3e_1,\\
&z_4e_3\cdot z_4e_3e_5 = z_4e_5\cdot z_4e_5e_3,\, z_4e_3\cdot z_4e_3e_4 = z_4e_4\cdot z_4e_4e_3.
\end{align*}

An easy simplification shows that neither new generators nor new relations are needed.

\medskip

From $z_5$ we have

\begin{align*}
z_5 &= (e_1,e_2,e_3,\overline e_3);(\overline e_1,e_4,e_5),(\overline e_2,\overline e_5,\overline e_4),\\
& z_5^1 = (e_4,e_5,e_2,e_3,\overline e_3);(\overline e_2,\overline e_5,\overline e_4),\\
& z_5^2 = (e_1,\overline e_5,\overline e_4,e_3,\overline e_3);(\overline e_1,e_4,e_5),\\
& z_5^4 = (e_1,e_2,e_3,\overline e_3);(\overline e_1,\overline e_2,\overline e_5,e_5),\\
& z_5^5 = (e_1,e_2,e_3,\overline e_3);(\overline e_1,e_4,\overline e_4,\overline e_2).
\end{align*}

Generators for $z_5$ are:
\begin{align*}
&z_5e_1,\, z_5e_2,\, z_5e_4,\, z_5e_5,\\
&z_5e_1e_4,\, z_5e_1e_5,\, z_5e_1e_2,\\
&z_5e_2e_1,\, z_5e_2e_5,\, z_5e_2e_4,\\
&z_5e_4e_1,\, z_5e_4e_2,\\
&z_5e_5e_1,\, z_5e_5e_2. \\
\end{align*}

We see $z_1^3=g_{5,1}^{1,3}\cdot z_5^1, \, z_1^3 = g_{5,1}^{2,3}\cdot z_5^2$, where
$$
\begin{array}{ll}
g_{5,1}^{1,3}=\left\{\begin{array}{lll}
e_2 & \mapsto & \overline e_5,\\
e_3 & \mapsto & \overline e_4,\\
e_4 & \mapsto & e_1,\\
e_5 & \mapsto & e_2;\\
\end{array}\right.
& 
g_{5,1}^{2,3}=\left\{\begin{array}{lll}
e_1 & \mapsto & e_1,\\
e_3 & \mapsto & \overline e_2,\\
e_4 & \mapsto & e_5,\\
e_5 & \mapsto & \overline e_4.
\end{array}\right.
\end{array}
$$

Relations for $z_5$ are:
\begin{align*}
&z_5e_1=1,\, z_5e_4=1,\, z_5e_5=1,\\
&z_5e_1e_4 = z_1e_3e_1,\, z_5e_1e_5 = z_1e_3e_2,\,z_5e_1e_2 = z_1e_2e_4 ,\\
&z_5e_2e_1 = z_1e_3e_1,\, z_5e_2e_5 = z_1e_3e_4,\, z_5e_2e_4 = z_1e_3e_5,\\
&z_5e_1\cdot z_5e_1e_4 = z_5e_4\cdot z_5e_4e_1,\,  z_5e_1\cdot z_5e_1e_5 = z_5e_5\cdot z_5e_5e_1,\\
& \quad\quad z_5e_1\cdot z_5e_1e_2 = z_5e_2\cdot z_5e_3e_1,\\
&z_5e_2\cdot z_5e_2e_5 = z_5e_5\cdot z_5e_5e_2,\, z_5e_2\cdot z_5e_2e_4 = z_5e_4\cdot z_5e_4e_2.
\end{align*}

An easy simplification shows that neither new generators nor new relations are needed.

\medskip

From $z_6$ we have

\begin{align*}
z_6 &= (e_1,e_2,e_3,\overline e_1);(\overline e_2,e_4,e_5),(\overline e_3,\overline e_5,\overline e_4),\\
& z_6^2 = (e_1,e_4,e_5,e_3,\overline e_1);(\overline e_3,\overline e_5,\overline e_4),\\
& z_6^3 = (e_1,e_2,\overline e_5,\overline e_4,\overline e_1);(\overline e_2,e_4,e_5),\\
& z_6^4 = (e_1,e_2,e_3,\overline e_1);(\overline e_2,\overline e_3,\overline e_5,e_5),\\
& z_6^5 = (e_1,e_2,e_3,\overline e_1);(\overline e_2,e_4,\overline e_4,\overline e_3).
\end{align*}

Generators for $z_6$ are:
\begin{align*}
&z_6e_2,\, z_6e_3,\, z_6e_4,\, z_6e_5,\\
&z_6e_2e_4,\, z_6e_2e_5,\, z_6e_2e_3,\\
&z_6e_3e_2,\, z_6e_3e_5,\, z_6e_3e_4,\\
&z_6e_4e_2,\, z_6e_4e_3,\\
&z_6e_5e_2,\, z_6e_5e_3. \\
\end{align*}

We see $z_2^4=g_{6,2}^{2,4}\cdot z_6^2, \, z_2^4 = g_{6,2}^{3,4}\cdot z_6^3$, where
$$
\begin{array}{ll}
g_{6,2}^{2,4}=\left\{\begin{array}{lll}
e_1 & \mapsto & e_1,\\
e_3 & \mapsto & e_5,\\
e_4 & \mapsto & e_2,\\
e_5 & \mapsto & e_3;\\
\end{array}\right.
& 
g_{6,2}^{3,4}=\left\{\begin{array}{lll}
e_1 & \mapsto & e_1,\\
e_2 & \mapsto & e_2,\\
e_4 & \mapsto & \overline e_5,\\
e_5 & \mapsto & \overline e_3.
\end{array}\right.
\end{array}
$$

Relations for $z_6$ are:
\begin{align*}
&z_6e_2, z_6e_4=1,\, z_6e_5=1,\\
&z_6e_2e_4 = z_2e_4e_2,\, z_6e_2e_5 = z_2e_4e_3,\, z_6e_2e_3 = z_2e_4e_5,\\
&z_6e_3e_2 = z_2e_4e_2,\, z_6e_3e_5 = z_2e_4e_3,\, z_6e_3e_4 = z_2e_4e_5,\\
&z_6e_2\cdot z_6e_2e_4 = z_6e_4\cdot z_6e_4e_2,\,  z_6e_2\cdot z_6e_2e_5 = z_6e_5\cdot z_6e_5e_2,\\
& \quad\quad z_6e_2\cdot z_6e_2e_3 = z_6e_3\cdot z_6e_3e_2,\\
&z_6e_3\cdot z_6e_3e_5 = z_6e_5\cdot z_6e_5e_3,\, z_6e_3\cdot z_6e_3e_4 = z_6e_4\cdot z_6e_4e_3.
\end{align*}

An easy simplification shows that neither new generators nor new relations are needed.

\medskip

We conclude 
\begin {align*}
\AM_{0,1,3} & = \left\langle z_1e_2e_4,\, z_1e_3e_5,\, z_1e_4e_5 \,\,\Big| \begin{array}{l}z_1e_4e_5=z_1e_2e_4\cdot z_1e_3e_5\cdot z_1e_2e_4,\\ 
z_1e_4e_5\cdot z_1e_3e_5=z_1e_2e_4\cdot z_1e_4e_5\end{array} \right\rangle\\
& = \langle z_1e_2e_4,\, z_1e_3e_5 \mid z_1e_3e_5\cdot z_1e_2e_4\cdot z_1e_3e_5=z_1e_2e_4\cdot z_1e_3e_5\cdot z_1e_2e_4 \rangle.
\end {align*}
\qed
\end{example}

\medskip

\noindent{\textbf{{Acknowledgments}}}

\medskip
\footnotesize

I am grateful to Warren Dicks for introducing me to algebraic mapping-class groups and to Luis Paris for his advices during a 
post-doc in Dijon.

\vskip -.5cm\null

\bibliographystyle{plain}
\bibliography{alg-present}

\end{document}